\documentclass[11pt]{article}
\usepackage{amssymb}
\usepackage{amsmath}

\begin{document}

\title{The geometry of product conjugate connections}

\markboth{{\small\it {\hspace{2cm} The geometry of product
conjugate connections}}}{\small\it{The geometry of product
conjugate connections \hspace{2cm}}}

\author{Adara M. Blaga and Mircea Crasmareanu}

\date{\small{\textit{Dedicated to Professor Vasile Cruceanu on the occasion of his 80th birthday}}\\
An. Stiint. Univ. Al. I. Cuza Iasi Math., 59(2013), no. 1, 73-84}

\maketitle

\begin{abstract}
Properties of pairs of product conjugate connections are stated
with a special view towards the integrability of the given almost
product structure. We define the analogous in product geometry of
the structural and the virtual tensors from the Hermitian geometry
and express the product conjugate connections in terms of these
tensors. Some examples from the geometry of a pair of
complementary distributions are discussed and for this case the
above structural and virtual tensors are expressed in terms of
O'Neill-Gray tensor fields.
\end{abstract}

\medskip

\noindent AMS Subject Classification: 53B05, 53C15.

\noindent Keywords and phrases: almost product structure;
(conjugate) linear connection; almost product Riemannian manifold;
structural and virtual tensor field.

\medskip

\section*{Introduction}

Fix $M$ a smooth, $n$-dimensional manifold for which we denote:
$C^{\infty }\left(M\right)$ -- the algebra of smooth real functions
on $M$, ${\mathfrak X}\left(M\right)$ -- the Lie algebra of vector
fields on $M$, $T_{s}^{r}\left(M\right)$ -- the $C^{\infty }\left(
M\right)$-module of tensor fields of $\left(r, s\right) $-type on
$M$. Usually $X, Y, Z, ...$ will be vector fields on $M$ and if
$T\rightarrow M$ is a vector bundle over $M$, then $\Gamma (T)$
denotes the $C^{\infty }$-module of sections of $T$ [e.g. $\Gamma
(TM)={\mathfrak X}(M)$].

\medskip

Let ${\mathcal C}(M)$ be the set of linear connections on $M$. Since
the difference of two linear connections is a tensor field of $(1,
2)$-type, it results that ${\mathcal C}(M)$ is a $C^{\infty
}(M)$-affine module associated to the $C^{\infty}(M)$-linear module
$T^1_2(M)$.

\medskip

Fix now $E$ an almost product structure on $M$, i.e. an endomorphism
of the tangent bundle such that $E^2=I_{{\mathfrak X}(M)}$. Then the
associated linear connections are:

\medskip

{\bf Definition 0.1} $\nabla \in {\mathcal C}(M)$ is an $E$-{\it
connection} if $E$ is covariant constant with respect to $\nabla $,
namely $\nabla E=0$. Let ${\mathcal C}_E(M)$ be the set of these
connections.

\medskip

In order to find the above set, let us consider after \cite[p.
342]{v:c} the maps
$$
\psi _E: {\mathcal C}(M)\rightarrow {\mathcal C}(M), \quad \chi
_E: T^1_2(M)\rightarrow T^1_2(M) \eqno(0.1)
$$
given by
$$
\psi _E(\nabla):=\frac{1}{2}\left(\nabla +E\circ \nabla \circ
E\right), \quad \chi _E(\tau):=\frac{1}{2}(\tau +E\circ \tau \circ
E). \eqno(0.2)
$$
So
$$
\left\{
\begin{array}{ll}
\psi _E(\nabla )_XY=\frac{\displaystyle 1}{\displaystyle 2}\left[\nabla _XY+E(\nabla _XEY)\right] \\
\chi _E(\tau)(X, Y)=\frac{\displaystyle 1}{\displaystyle
2}\left[\tau (X, Y)+E(\tau (X, EY))\right].
\end{array}
\right.  \eqno(0.3)
$$
Then, $\psi _E$ is a $C^{\infty }(M)$-projector on ${\mathcal C}(M)$
associated to the $C^{\infty }(M)$-linear projector $\chi _E$:
$$
\psi _E^2=\psi _E, \quad \chi _E^2=\chi _E, \quad \psi _E(\nabla
+\tau )=\psi _E(\nabla )+\chi _E(\tau ). \eqno(0.4)
$$

\smallskip

It follows that $\nabla E=0$ means $\psi _E(\nabla )=\nabla $ which
gives that ${\mathcal C}_E(M)=Im \psi _E$. This determines
completely ${\mathcal C}_E(M)$. Fix $\nabla _0$ arbitrary in
${\mathcal C}(M)$ and $\nabla $ in ${\mathcal C}_E(M)$. So, $\nabla
=\psi _E(\nabla ')$ with $\nabla '=\nabla _0+\tau $. In conclusion,
$\nabla =\psi _E(\nabla _0)+\chi _E(\tau )$; in other words,
${\mathcal C}_E(M)$ is the affine submodule of ${\mathcal C}(M)$
passing through the $E$-connection $\psi _E(\nabla _0)$ and having
the direction given by the linear submodule $Im \chi_E$ of
$T^1_2(M)$.

\medskip

Let us remark a decomposition (of arithmetic mean type) of it
\cite[p. 343]{v:c}:
$$
\psi _E(\nabla )=\frac{1}{2}(\nabla +C_E(\nabla )) \eqno(0.5)
$$
with {\it the conjugation map} $C_E: {\mathcal C}(M)\rightarrow {\mathcal C}(M)$:
$$
C_E(\nabla )_X:=E\circ \nabla _X\circ E. \eqno(0.6)
$$
Then {\it the product conjugate connection} $C_E(\nabla )$ measures
how far the connection $\nabla $ is from being an $E$-connection and
as it is pointed out in \cite[p. 343]{v:c}, $C_E$ is the affine
symmetry of the affine module ${\mathcal C}(M)$ with respect to the
affine submodule ${\mathcal C}_E(M)$, made parallel with the linear
submodule $\ker \chi _E$.

\medskip

The present paper is devoted to a large study of this new connection
$C_E(\nabla )$, since all above computations put in evidence its
r\^ole in the geometry of $E$. Therefore, the aim of our study is to
obtain several properties of it in both the general case and
Riemannian geometry. The first section is devoted to this
scope and after a general result connecting $\nabla $ and $C_E(\nabla )$, we treat two items: \\
i) the behavior of the product conjugate connections to a linear
change of almost product structure, \\
ii) the introduction in the product geometry of two tensor fields
previously considered in the complex geometry.\\
With respect to i) we arrive at two particular remarkable cases
concerning the recurrence of the given almost product structures,
while for ii) we derive some useful new identities.

\medskip

The second part of this paper is directed towards examples and to
put in our framework the geometry of (two complementary)
distributions. The most important case is when the considered
distributions are in a natural relationship with the initial linear
connection $\nabla $ or with the almost product structure $E$.

\medskip

In the last section we give some generalizations of the results
from the first part by adding an arbitrary tensor field of $(1,
2)$-type. All generalized product conjugate connections which form
a duality with the initial linear connection are determined.

\section{Properties of the product conjugate connection}

In what follows, for simplification we will denote by a superscript
$E$ the product conjugate connection of $\nabla$
$$
\nabla^{(E)}:=C_E(\nabla )=\nabla +E\circ \nabla E \eqno(1.1)
$$
and then
$$
\nabla ^{(E)}_XY=\nabla_ XY+E(\nabla _X EY-E(\nabla _XY))=E(\nabla _XEY). \eqno(1.2)
$$

The first properties of the product conjugate connection are stated in the
next proposition:

\medskip

{\bf Proposition 1.1} {\it Let $E$ be an almost product structure,
$\nabla $ a linear connection and $\nabla ^{(E)}$ the product
conjugate connection of $\nabla$. Then}:
\begin{enumerate}
  \item $\nabla ^{(E)}E=-\nabla E$; {\it it results that $\nabla \in {\mathcal C}_E(M)$ if and only if} $\nabla ^{(E)}\in {\mathcal C}_E(M)$.
  \item $\nabla $ {\it and $\nabla ^{(E)}$ are in duality}: $(\nabla
  ^{(E)})^{(E)}=\nabla $.
  \item $T_{\nabla ^{(E)}}=T_{\nabla }+E(d^{\nabla }E)$, {\it where $d^{\nabla }$ is the exterior covariant derivative induced by $\nabla $,
  namely $(d^{\nabla }E)(X, Y):=(\nabla _XE)Y-(\nabla _YE)X$; it results that for
  $\nabla \in {\mathcal C}_E(M)$, the connections $\nabla $ and $\nabla ^{(E)}$ have the same torsion}.
  \item $R_{\nabla ^{(E)}}(X, Y, Z)=E(R_{\nabla }(X, Y, EZ))$; {\it it results that
  $\nabla $ is flat if and only if $\nabla ^{(E)}$ is so}.
  \item {\it Assume that $(M, g, E)$ is an almost product Riemannian manifold, i.e. $g(EX, EY)=g(X,
  Y)$. Then $(\nabla ^{(E)}_Xg)(EY, EZ)=(\nabla _Xg)(Y, Z)$; it results
that $\nabla $ is a $g$-metric connection if and only if $\nabla
^{(E)}$ is so}.
\end{enumerate}

\medskip

{\bf Proof} 1. The main relations used here are
$$
\nabla ^{(E)}_XEY=E(\nabla _XY), \quad E(\nabla ^{(E)}_XY)=\nabla _XEY \eqno(1.3)
$$
and then
$$
(\nabla _XE)Y=\nabla _XEY-E(\nabla _XY)=E(\nabla ^{(E)}_XY)-\nabla ^{(E)}_XEY=-(\nabla ^{(E)}_XE)Y. \eqno(1.4)
$$
2. Although a direct proof can be provided by the formula $(0.6)$,
we prefer a proof here, in order to use $(1.1)$:
$$
(\nabla ^{(E)})^{(E)}=\nabla ^{(E)}+E\circ \nabla ^{(E)}E=\nabla
+E\circ \nabla E+E\circ (-\nabla E)=\nabla .
$$
3. A direct computation gives
$$
T_{\nabla ^{(E)}}(X, Y):=\nabla ^{(E)}_XY-\nabla
^{(E)}_YX-[X,Y]=E(\nabla _XEY)-E(\nabla _YEX)-[X,Y]=
$$
$$
=E(\nabla _XEY-\nabla _YEX)+T_{\nabla }(X,Y)-\nabla _XY+\nabla _YX:=
$$
$$
:=T_{\nabla }(X, Y)+E((\nabla _XE)Y-(\nabla _YE)X). \eqno(1.5)
$$
4.
$$
R_{\nabla ^{(E)}}(X, Y, Z):=\nabla ^{(E)}_X\nabla ^{(E)}_YZ-\nabla
^{(E)}_Y\nabla ^{(E)}_XZ-\nabla ^{(E)}_{[X,Y]}Z =
$$
$$
=\nabla ^{(E)}_XE(\nabla _YEZ)-\nabla ^{(E)}_YE(\nabla _XEZ)-E(\nabla _{[X,Y]}EZ)=
$$
$$
=E(\nabla _X\nabla _YEZ)-E(\nabla _Y\nabla _XEZ)-E(\nabla
_{[X,Y]}EZ):=E(R_{\nabla }(X, Y, EZ)). \eqno(1.6)
$$
5.
$$
(\nabla ^{(E)}_Xg)(V, W):=X(g(V, W))-g(\nabla ^{(E)}_XV, W)-g(V,
\nabla ^{(E)}_XW)=
$$
$$
=X(g(V, W))-g(E(\nabla _XEV), W)-g(V, E(\nabla _XEW))
$$
for any $X$, $V$ and $W\in {\mathfrak X}(M)$. With $V:=EY$ and
$W:=EZ$, we get
$$
(\nabla ^{(E)}_Xg)(EY, EZ)=X(g(EY, EZ))-g(E(\nabla _XY), EZ)-g(EY, E(\nabla _XZ))=
$$
$$
=X(g(Y, Z))-g(\nabla _XY, Z)-g(Y, \nabla _XZ):=(\nabla _Xg)(Y, Z).
\eqno(1.7)
$$
The above substitutions hold for $Y=EV$ and $Z=EW$. \quad $\Box$

\medskip

There are some direct consequences of these formulae: \\
i) if $\nabla $ is the Levi-Civita connection of $g$, then $\nabla ^{(E)}$ is also metric with respect to $g$, \\
ii) if $\nabla $ is the Levi-Civita connection of $g$ and in
addition $\nabla \in {\mathcal C}_E(M)$, then $\nabla ^{(E)}=\nabla
$ as the unique symmetric $g$-metric connection.

More generally, let $f\in Diff(M)$ be an automorphism of the
$G$-structure defined by $E$, i.e. $f_*\circ E=E\circ f_*$. If $f$
is an affine transformation for $\nabla $, namely $f_*(\nabla
_XY)=\nabla _{f_*X}f_*Y$, then $f$ is also affine transformation for
$\nabla ^{(E)}$.

\medskip

Two natural generalizations of the case $\nabla \in {\mathcal C}_E(M)$ are given by:

\medskip

{\bf Proposition 1.2} {\it Let $\nabla $ be a symmetric linear connection}. \\
i) {\it Assume that $E$ is $\nabla $-recurrent, i.e. $\nabla E=\eta
\otimes E$, where $\eta $ is a $1$-form. Then
$\nabla ^{(E)}$ is a semi-symmetric connection}. \\
ii) {\it Assume that $\nabla E=\eta \otimes I_{{\mathfrak X}(M)}$.
Then $\nabla ^{(E)}$ is a quarter-symmetric connection}.

\medskip

{\bf Proof} i) We have $\nabla ^{(E)}=\nabla +\eta \otimes I$ and from the item 3 of the
previous Proposition, we get $T_{\nabla ^{(E)}}=\eta \otimes I-I\otimes \eta $. \\
ii) It results that $\nabla ^{(E)}=\nabla +\eta \otimes E$ and, as
above, we get $T_{\nabla ^{(E)}}=\eta \otimes E-E\otimes \eta $.
\quad $\Box $

\medskip

The next subject consists of the behavior of $\nabla ^{(.)}$ for
families of almost product structures. Let $E_1$ and $E_2$ be two
almost product structures and consider the pencil of $(1, 1)$-tensor
fields $E_{\alpha , \beta }:=\alpha E_1+\beta E_2$ with
$\alpha $ and $\beta \in \mathbb{R}$. In order that $E_{\alpha , \beta}$ to be an almost product structure there are necessary two conditions: \\
1) $E_1$ and $E_2$ be skew-commuting structures: $E_1E_2=-E_2E_1$, \\
2) $(\alpha , \beta )$ belongs to the unit circle $S^1$: $\alpha^2+\beta^2=1$. \\
Then:
$$
\nabla ^{(E_{\alpha , \beta })}_XY=\alpha ^2\nabla ^{(E_1)}_XY+\beta ^2\nabla ^{(E_2)}_XY+\alpha \beta [E_1(\nabla _XE_2Y)+E_2(\nabla _XE_1Y)] \eqno(1.8)
$$
and there can be distinguished two remarkable particular cases: \\
i) if $E_1$ and $E_2$ are recurrent with respect to $\nabla $ with
the same $1$-form of recurrence: $\nabla E_i=\eta \otimes E_i$, then
the product conjugate connections coincide $\nabla ^{(E_1)}\equiv
\nabla ^{(E_2)}=:\nabla^{(E_{12})}$ and it follows the invariance of
$\nabla ^{(E)}$:
$$
\nabla ^{(E_{\alpha , \beta})}=\nabla ^{(E_{12})}, \eqno(1.9)
$$
ii) assume that the triple $(\nabla , E_1, E_2)$ is a
mixed-recurrent structure:  $\nabla E_i=\eta \otimes E_j$ with
$i\neq j$. Then $\nabla $ is the average of the two product
conjugate connections, $\nabla=\frac{\displaystyle 1}{\displaystyle
2}(\nabla ^{(E_1)}+\nabla ^{(E_2)})$ and
$$
\nabla ^{(E_{\alpha , \beta})}=\nabla +(\alpha ^2-\beta ^2)\eta \otimes E_1E_2. \eqno(1.10)
$$

\smallskip

The last subject of this section treats two tensor fields associated
to an almost product structure. The paper \cite{v:k} introduces the
structural and virtual tensor fields of an almost complex structure.
Turning into our framework, let us consider for
a pair $(\nabla , E)$ the tensor fields of $(1, 2)$-type: \\
1) {\it the structural tensor field}
$$
C^E_{\nabla }(X,Y):=\frac{1}{2}[(\nabla _{EX}E)Y+(\nabla _{X}E)EY]
\eqno(1.11)
$$
2) {\it the virtual tensor field}
$$
B^E_{\nabla }(X, Y):=\frac{1}{2}[(\nabla _{EX}E)Y-(\nabla
_{X}E)EY]. \eqno(1.12)
$$
From the item 1 of the first Proposition it results that both
these tensor fields are skew-symmetric with respect to the product
conjugation of connections:
$$
C^{E}_{\nabla ^{(E)}}=-C^E_{\nabla }, \ \ B^{E}_{\nabla ^{(E)}
}=-B^E_{\nabla }. \eqno(1.13)
$$
Also
$$
C^E_{\nabla }(EX, EY)=C^E_{\nabla }(X, Y), \ B^E_{\nabla }(EX,
EY)=-B^E_{\nabla }(X, Y). \eqno(1.14)
$$

\smallskip

The importance of these tensor fields for our study is given by
the following straightforward relation:
$$
\nabla ^{(E)}=\nabla -C^{E}_{\nabla }+B^{E}_{\nabla }. \eqno(1.15)
$$

\smallskip

Recall after \cite{c:mz} that two linear connections are called
{\it projectively equivalent} if there exists a $1$-form $\tau $
such that:
$$
\nabla '=\nabla +\tau \otimes I+I\otimes \tau . \eqno(1.16)
$$
A straightforward calculus gives that $C^E$ is invariant for
projectively changes $(1.16)$ while for $B^E$ we have:
$$
(B^E_{\nabla '}-B^E_{\nabla })(X, Y)=\tau (EY)EX-\tau (Y)X.
\eqno(1.17)
$$

\section{Invariant distributions}

Let $\mathcal{D}\subset TM$ be a fixed distribution considered as
a vector subbundle of $TM$.

\medskip

{\bf Definition 2.1} i) $\mathcal{D}$ is called $E$-{\it
invariant} if $X\in \Gamma (\mathcal{D})$ implies $EX \in \Gamma (\mathcal{D})$. \\
ii)(\cite[p. 210]{b:l}) The linear connection $\nabla $ {\it
restricts to} $\mathcal{D}$ if $Y\in \Gamma (\mathcal{D})$ implies
$\nabla _XY \in \Gamma (\mathcal{D})$, for any $X\in \Gamma (TM)$.

\medskip

If ${\nabla }$ restricts to ${\mathcal D}$, then $\nabla $ may be
considered as a connection in the vector bundle ${\mathcal D}$. From
this fact, in \cite[p. 7]{b:f} a connection which restricts to
${\mathcal D}$ is called {\it adapted to} ${\mathcal D}$.

\medskip

{\bf Proposition 2.2} {\it If the distribution $\mathcal{D}$ is
$E$-invariant and the linear connection $\nabla $ restricts to
$\mathcal{D}$, then $\nabla ^{(E)}$ also restricts to}
$\mathcal{D}$.

\medskip

{\bf Proof} Fix $Y\in \Gamma (\mathcal{D})$. Then $EY\in \Gamma
(\mathcal{D})$ and for any $X\in \Gamma( TM)$ we have $\nabla _XY
\in \Gamma (\mathcal{D})$. Therefore, $\nabla ^{(E)}_XY=E(\nabla
_XEY)\in \Gamma (\mathcal{D})$. \quad $\Box$

\medskip

A more general notion like restricting to a distribution is that of
geodesically invariance \cite[p. 118]{b:l}. The distribution
$\mathcal{D}$ is $\nabla $-\textit{geodesically invariant} if for
every geodesic $\gamma :[a, b]\rightarrow M$ of $\nabla $ with
$\dot{\gamma }(a)\in \mathcal{D}_{\gamma (a)}$ it follows
$\dot{\gamma }(t)\in \mathcal{D}_{\gamma (t)}$ for any $t\in [a,
b]$. The cited book gives a necessary and sufficient condition for a
distribution $\mathcal{D}$ to be $\nabla $-geodesically invariant:
for any $X$ and $Y\in \Gamma (\mathcal{D})$, {\it the symmetric
product} $\langle X:Y\rangle:=\nabla _XY+\nabla _YX$ to belong to
$\Gamma (\mathcal{D})$ or equivalently, for any $X\in \Gamma
(\mathcal{D})$ to have $\nabla _XX\in \Gamma (\mathcal{D})$.

\medskip

The following result is a direct consequence of definitions:

\medskip

{\bf Proposition 2.3} {\it If the distribution $\mathcal{D}$ is
$E$-invariant and the linear connection $\nabla $ restricts to
$\mathcal{D}$, then $\mathcal{D}$ is geodesically invariant for}
$\nabla ^{(E)}$.

\medskip

{\bf Example 2.4} Assume that the tangent bundle $TM$ admits a
decomposition
$$
TM=V\oplus H \eqno(2.1)
$$
into {\it vertical} and {\it horizontal vectors}. Let ${\mathcal
D}_v=\Gamma (V)$ respectively, ${\mathcal D}_h=\Gamma (H)$ and the
corresponding projectors $v$ and $h$. Then $E=h-v$ is an almost
product structure and both ${\mathcal D}_*$ are $E$-invariant. As it
is proved in \cite{b:f}, the almost product structures are in a
natural relationship with decompositions of $(2.1)$-type and so,
$E=h-v$ is the prototype of all possible almost product structures.

\medskip

The product conjugate connection of $\nabla $ is
$$
\nabla ^{(E)}_XY=h(\nabla _XhY)-h(\nabla _XvY)-v(\nabla
_XhY)+v(\nabla _XvY) \eqno(2.2)
$$
and then we have:

\medskip

{\bf Proposition 2.5} {\it If $\nabla ^{(E)}$ is torsion-free, then
$E$ is integrable, which means that ${\mathcal D}_h$ and ${\mathcal
D}_v$ are involutive distributions}.

\medskip

{\bf Proof} From $(2.2)$ we get
$$
h[X, Y]+v[X, Y]=\nabla ^{(E)}_XY-\nabla ^{(E)}_YX=h(\nabla
_XhY-\nabla _YhX)+v(\nabla _XvY-\nabla _YvX)
$$
and then
$$
h[X, Y]=\nabla _XhY-\nabla _YhX, \quad v[X, Y]=\nabla _XvY-\nabla
_YvX.
$$
With $X\rightarrow vX$ and $Y\rightarrow vY$ in the first relation
above it follows $h[vX, vY]=0$ and the change $X\rightarrow hX$ and
$Y\rightarrow vY$ in the second relation yields $v[hX, hY]=0$. \quad
$\Box $

\medskip

We have: \\
1) $\nabla $ restricts to ${\mathcal D}_h$ means $v(\nabla
_XhY)=0$ and $h(\nabla _XhY)=\nabla _XhY$, \\
2) $\nabla $ restricts to ${\mathcal D}_v$ means $h(\nabla _XvY)=0$
and $v(\nabla _XvY)=\nabla _XvY$.

A straightforward computation gives that the general $\nabla ^{(E)}$
of $(2.2)$ restricts to ${\mathcal D}_h$ and ${\mathcal D }_v$.
Moreover, if $\nabla $ restricts to both ${\mathcal D}_*$, then
$$
\nabla ^{(E)}_XY=\nabla _XhY+\nabla _XvY=\nabla _XY \eqno(2.3)
$$
and so $\nabla \in {\mathcal C}_E(M)$. Let us remark that the above
connection $(2.3)$ is exactly {\it the Schouten connection} of the
pair $(h, v)$ \cite[p. 10]{f:ip}:
$$
\nabla _XY=h(\nabla _XhY)+v(\nabla _XvY). \eqno(2.4)
$$

\smallskip

{\bf Example 2.5} Let $M$ be a vector bundle over the manifold $U$
by $\pi :M\rightarrow U$ and $V=\ker T\pi $ {\it the vertical
bundle} over $M$. Let also $H_1$ and $H_2$ be two {\it horizontal
bundles} in the decomposition $(2.1)$ and denote $h_1$ respectively,
$h_2$ their projectors. Then the skew-symmetry $E_1E_2=-E_2E_1$ for
the corresponding almost product structures of the previous example
means the skew-symmetry $h_1h_2=-h_2h_1$.

\medskip

Sometimes, a complementary distribution $H$ to the above vertical
subbundle $V$ is called {\it Ehresmann connection}, while if $M$ is
exactly the tangent bundle $TU$, then $H$ is called {\it nonlinear
connection} \cite{b:fa}.

\medskip

{\bf Example 2.6} Let $E$ be the almost product structure of example
2.4. Our next step is to express the Kirichenko tensor fields in
terms of the projectors $h$, $v$:

\medskip

{\bf Proposition 2.7} {\it The structural and virtual tensor
fields of $E=h-v$ are}:
$$
\left\{
  \begin{array}{ll}
    C^{h-v}_{\nabla }(X,Y)=2[h(\nabla _{vX}vY)+v(\nabla _{hX}hY)] \\
    B^{h-v}_{\nabla }(X,Y)=-2[h(\nabla _{hX}vY)+v(\nabla _{vX}hY)].
  \end{array}
\right. \eqno(2.5)
$$

\smallskip

{\bf Proof} From $(1.14)$ we get
$$
\left\{
  \begin{array}{ll}
    C^{h-v}_{\nabla }(hX, vY)=-C^{h-v}_{\nabla }(vX, hY) \\
    B^{h-v}_{\nabla }(hX, hY)=-B^{h-v}_{\nabla }(vX, vY).
  \end{array}
\right.  \eqno(2.6)
$$
By making $X\rightarrow vX$ in the first relation and $X\rightarrow
hX$ in the second one, it results
$$
\left\{
  \begin{array}{ll}
    C^{h-v}_{\nabla }(hX, vY)=0=C^{h-v}_{\nabla }(vX, hY) \\
    B^{h-v}_{\nabla }(hX, hY)=0=B^{h-v}_{\nabla }(vX, vY)
  \end{array}
\right. \eqno(2.7)
$$
and then
$$
\left\{
  \begin{array}{ll}
    C^{h-v}_{\nabla }(X, Y)=C^{h-v}_{\nabla }(hX, hY)+C^{h-v}_{\nabla }(vX, vY) \\
    B^{h-v}_{\nabla }(X, Y)=B^{h-v}_{\nabla }(hX, hY)+B^{h-v}_{\nabla }(vX, vY).
  \end{array}
\right.  \eqno(2.8)
$$
But
$$
\left\{
  \begin{array}{ll}
    C^{h-v}_{\nabla }(hX, hY)=2v(\nabla _{hX}hY) \\
    C^{h-v}_{\nabla }(vX, vY)=2h(\nabla _{vX}vY)
  \end{array}
\right. \eqno(2.9)
$$
and
$$
\left\{
  \begin{array}{ll}
    B^{h-v}_{\nabla }(hX, vY)=-2h(\nabla _{hX}vY) \\
    B^{h-v}_{\nabla }(vX, hY)=-2v(\nabla _{vX}hY)
  \end{array}
\right.  \eqno(2.10)
$$
and then we have the conclusion. \quad $\Box $

\medskip

Let us recall the well-known {\it fundamental tensor fields} of O'Neill-Gray:
$$
\left\{
  \begin{array}{ll}
    T(X, Y)=h(\nabla _{vX}vY)+v(\nabla _{vX}hY) \\
    A(X, Y)=v(\nabla _{hX}hY)+h(\nabla _{hX}vY).
  \end{array}
\right. \eqno(2.11)
$$
Then, a comparison of last two equations yields
$$
\left\{
  \begin{array}{ll}
    C^{h-v}_{\nabla }(X, Y)=2[T(X, vY)+A(X, hY)] \\
    B^{h-v}_{\nabla }(X, Y)=-2[T(X, hY)+A(X, vY)]
  \end{array}
\right. \eqno(2.12)
$$
a fact which justifies the second name of $T$ and $A$ as {\it
invariants} of $(2.1)$ \cite[p. 9]{f:ip}.

\section{Generalized product conjugate connections}

In this section we present a natural generalization of the product
conjugate connection.

\medskip

{\bf Definition 3.1} A {\it generalized product conjugate
connection} of $\nabla $ is
$$
\nabla ^{(E, C)}=\nabla ^{(E)}+C \eqno(3.1)
$$
with $C\in T^1_2(M)$ arbitrary.

\medskip

Since the duality $\nabla \leftrightarrow \nabla ^{(E)}$ is a main
feature of $\nabla ^{(E)}$, let us search for tensor fields $C$ such
that $(\nabla ^{(E, C)})^{(E, C)}=\nabla $. From
$$
(\nabla ^{(E, C)})^{(E, C)}_XY=\nabla _XY+E(C(X, EY))+C(X, Y)
\eqno(3.2)
$$
it results that we are interested in finding solutions $C$ to
$$
E(C(X, EY))+C(X, Y)=0. \eqno(3.3)
$$

Let us remark that: \\
i) $C_0=\nabla E$ is a particular solution of $(3.3)$, \\
ii) if $C$ is a solution, then $E\circ C$ is also a solution.

So, let us search the duality property for
$$
\nabla ^{(E, \lambda , \mu )}=\nabla ^{(E)}+\lambda \nabla E +\mu
E(\nabla E)=(1+\mu )\nabla ^{(E)}+\lambda \nabla E \eqno(3.4)
$$
with $\lambda$ and $\mu \in \mathbb{R}$.

\medskip

{\bf Proposition 3.2} {\it The duality $\nabla \leftrightarrow
\nabla ^{(E, \lambda , \mu )}$ holds only for the pairs $(\lambda
, \mu )\in \{(0, 0), (0, -2), (1, -1), (-1, -1)\}$}.

\medskip

{\bf Proof} From
$$
(\nabla ^{(E, \lambda , \mu })^{(E, \lambda , \mu )}_XY=[(1+\mu
)^2+\lambda ^2]\nabla _XY+2\lambda (1+\mu )E(\nabla _XY)
$$
it results the system
$$
  \begin{cases}
    (1+\mu )^2+\lambda ^2=1 \\
    \lambda (1+\mu )=0.
  \end{cases}
$$
which has the above solutions.

Let us point out that
$$
  \begin{cases}
    \nabla ^{(E, 0, 0)}=\nabla ^{(E)} \\
    \nabla ^{(E, 0, -2)}=-\nabla ^{(E)}  \\
    \nabla ^{(E, 1, -1)}=\nabla \\
    \nabla ^{(E, -1, -1)}=-\nabla
  \end{cases}
$$
which confirm our result. \quad $\Box $

\medskip

Returning to the general case $(3.1)$, let us present the
generalizations of some relations from Proposition 1.1:\\
1. $\nabla ^{(E, C)}E=-\nabla E+C(\cdot , E\cdot )-E\circ C$. Then
$\nabla \in {\mathcal C}_E(M)$ if and only if
$\nabla ^{(E, \lambda \nabla E+\mu E\circ \nabla E)}\in {\mathcal C}_E(M)$ with $\lambda $ and $\mu $ arbitrary real numbers.\\
2. the discussion above. \\
3. $T_{\nabla ^{(E, C)}}=T_{\nabla }+E(d^{\nabla }E)+2C_{skew}$,
where $C_{skew}$ is the skew-symmetric part of $C$, i.e.
$2C_{skew}(X, Y)=C(X, Y)-C(Y, X)$. So, if $C$ is symmetric and
$\nabla \in {\mathcal C}_E(M)$,
then $\nabla $ and $\nabla ^{(E, C)}$ have the same torsion.\\
4. $R_{\nabla ^{(E, C)}}(X, Y)Z=E(R_{\nabla }(X, Y)EZ)+C(X,
E(\nabla _YEZ))-C(Y, E(\nabla _XEZ))-C([X, Y], Z)+E(\nabla _YE(
C(Y, Z)))-E(\nabla _YE(C(X, Z)))$.

\medskip

{\bf Acknowledgement.} The first author acknowledges the support by the research
grant PN-II-ID-PCE-2011-3-0921.

\begin{tabbing}
Department of Mathematics and Computer Science \\
West University of Timi\c{s}oara \\
Bld. V. P\^{a}rvan nr. 4, 300223  Timi\c{s}oara \\
Rom\^{a}nia \\
adara@math.uvt.ro\\
\end{tabbing}

\begin{tabbing}
Faculty of Mathematics \\
University "Al. I. Cuza" \\
Ia\c si, 700506 \\
Rom\^ania \\
mcrasm@uaic.ro \\

\smallskip

\noindent http://www.math.uaic.ro/$\sim$mcrasm
\end{tabbing}


\begin{thebibliography}{99}

\bibitem{b:f} A. Bejancu, H. R. Farran, {\it Foliations and geometric structures},
Mathematics and Its Applications 580, Springer, Dordrecht, 2006. MR2190039 (2006j:53034)

\bibitem{b:fa} A. Bejancu, H. R. Farran, {\it Riemannian Metrics on the Tangent
Bundle of a Finsler Submanifold}, An. \c Stiin\c t. Univ. Al. I.
Cuza Ia\c si. Mat., 56(2010), no. 2, 429-436. DOI
10.2478/v10157-010-0030-8

\bibitem{b:l} F. Bullo, A. D. Lewis, {\it Geometric control of mechanical systems.
Modeling, analysis, and design for simple mechanical control systems}, Texts in Applied Mathematics, 49, Springer-Verlag, New York, 2005. MR2099139 (2005h:70030)

\bibitem{c:mz} O. Calin, H. Matsuzoe, J. Zhang, {\it Generalizations of conjugate
connections}, Proceedings of 9th International Workshop on Complex
Structures, Integrability, and Vector Fields, 2009, 26-34.

\bibitem{v:c} V. Cruceanu, {\it Almost hyperproduct structures on manifolds}, An. \c Stiin\c t. Univ. Al. I.
Cuza Ia\c si Mat., 48(2002), no. 2, 337-354. MR2007448 (2004h:53040)

\bibitem{f:ip} M. Falcitelli, S. Ianu\c s, A. M. Pastore, {\it Riemannian submersions and related topics},
 World Scientific Publishing Co., Inc., River Edge, NJ, 2004. MR2110043 (2005k:53036)

\bibitem{v:k} V. F. Kirichenko, {\it Method of generalized Hermitian
geometry in the theory of almost contact manifold}, Itogi Nauki i
Tekhniki, Problems of geometry, 18(1986), 25-71; translated in J.
Soviet. Math., 42(1988), no.5, 1885-1919. MR0895367 (88g:53042)

\end{thebibliography}
\end{document}